\documentclass[10pt]{amsart}
\usepackage{amssymb,amsmath, amscd,amsthm}
\usepackage[all]{xy}

\newcommand{\map}[2]{\,{:}\,#1\!\longrightarrow \! #2}

\newcommand{\Pro}{\mathbf{P}}

\newcommand{\pr}{\mathrm{pr}}
\newcommand{\tw}{\mathrm{tw}}

\newcommand{\Aff}{\mathbf{A}}

\newcommand{\sm}{\mathrm{sm}}

\newcommand{\Sp}{\mathbf{Sp}}
\newcommand{\PSp}{\mathbf{PSp}}
\newcommand{\C}{\mathbf{C}}

\newcommand{\Jac}{\mathrm{Jac}}
\newcommand{\curl}[1]{\mathcal{#1}}
\newcommand{\Spec}{\mathrm{Spec}}
\newcommand{\Gal}{\mathrm{Gal}}

\newcommand{\cH}{\mathrm{H}}

\newcommand{\tor}{\mathfrak{t}}

\newcommand{\pair}[2]{ \left[ \begin{array}{cc} #1 \\ #2 \\ \end{array} \right]}

\newcommand{\Q}{\mathbf{Q}}
\newcommand{\an}{\mathrm{an}}
\newcommand{\hol}{\mathrm{hol}}

\newcommand{\Disc}{\mathrm{Discr}}
\newcommand{\Res}{\mathrm{Res}}
\newcommand{\Sym}{\mathrm{Sym}}

\newcommand{\Proj}{\mathrm{Proj}}

\newcommand{\R}{\mathbf{R}}

\newcommand{\Z}{\mathbf{Z}}

\DeclareMathOperator{\Hom}{Hom}

\newcounter{thm}

\theoremstyle{plain}

\newtheorem{Theorem}[thm]{Theorem}
\newtheorem{Proposition}[thm]{Proposition}

\theoremstyle{definition}
\newtheorem{Definition}[thm]{Definition}

\theoremstyle{remark}

\begin{document}
\title{Torelli Revisited}
\author{Stephen Meagher}
\email{sjmeagher@gmail.com}

\begin{abstract}
Assume given an Abelian variety which is a geometric Jacobian.
We give an invariant of the moduli point of the Abelian variety which determines the minimal extension of its ground field over which it is a Jacobian.
For genus $3$, we express this invariant in terms of weight 18 Siegel and Teichm\"uller modular forms and periods of the Abelian variety. These results answer questions of Serre posed in letters to Top from 2003 about Abelian $3$-folds with indecomposable principal polarisations. The results we obtain hold in any characteristic other than $2$ and thus generalise previous work of Ritzenthaler, Lachaud and Zykin on the questions in Serre's letters.
\end{abstract}

\maketitle

\section{Introduction}
Let $(A,a)$ be an Abelian $g$-fold with an indecomposable principal polarisation $a$ over a field $k$; let $k^s$ be the separable closure of $k$. Assume that there is a smooth irreducible genus $g$ curve $C$ over $k^s$ whose Jacobian variety satisfies
\[ (\Jac(C), \lambda_\Theta) \cong (A,a) \otimes k^s, \]
where $\lambda_\Theta$ is the polarisation of $\Jac(C)$ induced by the theta divisor. (For example if $g \leq 3$ then such a $C$ always exists \cite{oortUeno}).

Now assume $k$ is not of characteristic $2$. For any character $\epsilon \map{\Gal(k^s/k)}{\{ \pm 1 \}}$ we let $(A,a)_\epsilon$ be the indecomposably principally polarised Abelian variety obtained by twisting $(A,a)$ over $(k^s)^{\ker(\epsilon)}$ by the $[-1]$ automorphism. Thus
if $\ker(\epsilon) = \Gal(k^s/k)$ then $(A,a)_\epsilon = (A,a)$. Serre has observed (\cite{Lauter} Th\'eor\`emes 4 \& 5) the following refinement of the Torelli theorem: there is a unique character $\epsilon \map{\Gal(k^s/k)}{\{ \pm 1\}}$ and a unique $C/k$ so that
\[ (\Jac(C), \lambda_\Theta) \cong_k (A,a)_{\epsilon}. \]
Moreover $\epsilon$ is trivial if $C$ is hyperelliptic. Thus there is a $D_{(A,a)} \in k$ such that $k(\sqrt{D_{(A,a)}}) = (k^s)^{\ker(\epsilon)}$ is the minimal field extension of $k$ over which $(A,a)$ is a Jacobian. This motivates
\begin{Definition} \label{twistDef} Let $R$ be a $\Z[1/2]$ algebra.
Let $(A,a)/\Spec(R)$ be a family\footnote{i.e. a commutative group scheme which is smooth and proper over $\Spec(R)$ and whose geometric fibres are connected Abelian $g$-folds.} of principally polarised Abelian $g$-folds whose geometric fibres are Jacobians. A function $\Delta \in R$ is called a twisting function if for all fields $k$ and all points $x \in \Spec(R)(k)$ we have: \\
i) $\Delta(x) = 0$ if and only if $(A,a)_x$ is the Jacobian of a hyperelliptic curve over $k$ \\
ii) $(A,a)_x$ is a Jacobian if and only if $\Delta(x)$ is a square in $k$. 
\end{Definition}

Our main result is

\begin{Theorem} \label{twist}
Let $S$ be a scheme over which $2$ is invertible. Let $(A,a)/S$ be a family of principally polarised Abelian $g$-folds whose geometric fibres are Jacobians. If either $g = 3$ or the geometric fibres of $(A,a)/S$ are all non-hyperelliptic Jacobians, then there exists an affine open cover $\{U_i \}$ of $S$ and twisting functions $\Delta_i \in \cH^0(U_i, \mathcal{O}_{U_i})$ for $(A,a)_{U_i}$. 
\end{Theorem}

\textbf{Remark:} For $g > 3$ the hyperelliptic locus has codimension greater than $1$ in the moduli space of curves and is therefore not defined locally by a single equation, so in some sense Theorem \ref{twist} is sharp. However interesting relations may hold between local equations for the hyperelliptic locus and the twisting functions, but we are ignorant on this point. \\

Let $\curl{M}_3$ be the moduli stack of smooth genus three curves and let $\curl{A}_{3,1}^i$ be the moduli
stack of Abelian threefolds equipped with an indecomposable principal polarisation. The Torelli morphism 
\[ \tor \map{\curl{M}_3}{\curl{A}_{3,1}^i} \]
associates to a family\footnote{i.e. $C/S$ is a smooth, proper morphism whose geometric fibres are irreducible genus $3$ curves.} of genus three curves $C/S$ its family of principally polarised Jacobians $(\Jac(C)/S, \lambda_\Theta)$. The classical Torelli theorem implies that this map is a surjection on $k^s$ valued points. The refined Torelli theorem states that the fibre of $\tor$ above an Abelian
threefold $(A,a)$ defined over $k$ has splitting field given by $(k^s)^{\ker(\epsilon)}$ for the character $\epsilon \map{\Gal(k^s/k)}{\mu_2}$
as above. In particular Theorem 2 means that the Torelli morphism is ``$2$ to $1$ and ramified along the hyperelliptic locus over $\Z[1/2]$'', although we do not define this notion for stacks.

\begin{Theorem} \label{chidisc} Let $\chi_{18}$ be the modular form obtained by the product of the $36$ even theta nulls. Then $\chi_{18}$ extends to a modular form on $\mathcal{A}_{3,1}^i$ with coefficients in $\Z$. Moreover $\chi_{18}$ vanishes on the hyperelliptic locus of $\curl{A}_{3,1}^i$. Let
	 $\Disc$ be the weight $9$ Katz-Teichm\"uller modular form defined in Section \ref{disc}. Then
	 $\Disc$ vanishes on the hyperelliptic locus of $\curl{M}_3$ and $\tor^*\chi_{18} = \gamma
	 \Disc^2$ for some $\gamma \in \Z^*$.
\end{Theorem}

Combining Theorems \ref{twist} and \ref{chidisc} we obtain
\begin{Theorem} \label{twistchi}
	Let $R \subset \C$ be a ring containing $1/2$ and let $(A,a)$ be a principally polarised Abelian threefold over $R$. Assume that $\Omega_{A/R}$ is free and that there is a twisting function $\Delta \in R$ for $(A,a)$. Let $\{ \eta_1, \dots, \eta_6 \}$ be a basis for $\cH_1(A_\C^\an, \Z)$ and let $\xi_1, \xi_2, \xi_3$ be a basis for $\Omega_{A/R}(R)$. Let $\Omega_1$ be the matrix of integrals of the $\xi_i$ along the $\eta_j$ with $j \in \{1, 2, 3\}$ and let $\Omega_2$ be the matrix of integrals of the $\xi_i$ along the $\eta_j$ with $j \in \{4,5,6\}$.	Let $\tau = \Omega_2 \cdot \Omega_1^{-1}$ and let $\chi_{18}^{\hol}(\tau)$ denote the product of the $36$
	even theta nulls evaluated at $\tau$. Let $\gamma$ be as in Theorem \ref{chidisc}. Then for each maximal ideal $\mathfrak{m} \subset R$
	\[ (2\pi)^{54} \frac{\chi_{18}^{\hol}(\tau)}{( \det(\Omega_1))^{18}} + \mathfrak{m} = -\gamma \Delta(A,a) + \mathfrak{m}, \]
	up to a square unit in $R/\mathfrak{m}$.
\end{Theorem}

\textbf{Remarks: 1)} The hypotheses of Theorem \ref{twistchi} hold automatically for $R \subset \C$ a local $\Z[1/2]$ algebra. \\ 
\textbf{2)} The modular forms $\chi_{18}$ and $\Disc$ have long histories going back to Klein \cite{Klein} and Igusa \cite{igusa}. Lachaud, Ritzenthaler and Zykin \cite{ritz, ritz2} have also given proofs of Theorems \ref{chidisc} and \ref{twistchi} in characteristic $0$ without reference to twisting functions. They calculate explicitly that $\gamma = -1$. The main novelty of our treatment is a proof of Theorem \ref{twist} in genus $3$ which does not refer to modular forms. We give a second proof of Theorem \ref{twist} for $g=3$ as a corollary of Theorem \ref{chidisc}. Theorems \ref{chidisc} and \ref{twistchi} answer questions of Serre in letters to Top \cite{serreLet}. Explicit formulae for special types of twisting functions have been given in \cite{AuerTop} and \cite{HLP}. A suggestion about candidate modulur forms for analogues of Theorems \ref{chidisc} and \ref{twistchi} in genus 4 is given in the final section of \cite{ritz2}.

\section{Moduli spaces}

	Let $\curl{M}_g$ be the category fibred in groupoids over schemes whose objects are smooth proper morphisms of relative dimension $1$ with smooth connected genus $g$ curves for geometric fibres. Let $\curl{A}_{g,1}$ be the category fibred in groupoids over schemes whose objects are Abelian schemes of relative dimension $g$ equipped with a principal polarisation. The Torelli morphism is the morphism given by
	\[ \tor \map{\curl{M}_{g}}{\curl{A}_{g,1}} : C/S \mapsto (\Jac(C)/S, \lambda_\Theta). \]
	The categories $\curl{M}_g$ and $\curl{A}_{g,1}$ are Deligne-Mumford stacks (\cite{DelMum, FaltChai}) and the Torelli morphism is a morphism of fibred categories. \\

	For each $N \in \Z$ we fix a primitive $N$th root of unity $\zeta_N \in \bar{\Q}$ such that if $M$ and $N$ are coprime integers, $\zeta_{MN} = \zeta_M \zeta_N$. We consider only  symplectic level $N$ structures which send the determinant of a standard basis to $\zeta_N$. Let $\curl{M}_{g,N}$ be the functor from $\Z[1/2N, \zeta_N]$-schemes to sets given by
	\begin{eqnarray*}
		S & \mapsto & \{ (C, \alpha) \mid C \in \curl{M}_g(S) \text{ and } \alpha : \Jac(C)[N] \cong (\Z/N)^g_S \times (\Z/N)^g_S \text{ is symp. } \}/_{\cong_S}.
	\end{eqnarray*}	
	If $N \geq 3$ then $\curl{M}_{g,N}$ is represented by a smooth irreducible scheme over $\Z[1/2N, \zeta_N]$ (\cite{oortSteen} Theorem 1.8, \cite{Popp} p134 Theorem 10.10, p142 Remark 2, p104 Theorem 8.11; for smoothness \cite{Harris} p103 Lemma 3.35; for irreducibility \cite{DelMum}).
	
	Likewise for $N \geq 3$ the functor from $\Z[1/2N, \zeta_N]$-schemes to sets given by
	\begin{eqnarray*}
		S & \mapsto & \{ (A,a, \alpha) \mid (A,a) \in \curl{A}_{g,1}(S) \text{ and }\alpha : A[N] \cong (\Z/N)^g_S \times (\Z_N)^g_S \text{ is symp. } \}/_{\cong_S},
	\end{eqnarray*}
	is represented by an irreducible smooth projective scheme $\curl{A}_{g,1,N}$ over $\Z[ 1/2N, \zeta_N]$ (\cite{GIT} p139 Theorem 7.9; for smoothness see \cite{OortNordic} p242 Theorem 2.33; for irreducibility see \cite{FaltChai}). 

For lack of reference we prove
\begin{Proposition} \label{torsor}
Let $(A,a)/S$ be a principally polarised Abelian scheme of relative dimension $g$. Let $N \geq 1$ be an integer which is invertible over $S$ and let $\textbf{S}(A,a,N)$ be the the functor from $S$-schemes to sets whose $T$-valued points are given by
\[  \textbf{S}(A,a,N)(T) = \{ \alpha : A_T[N] \cong (\Z/N)^{2g}_T \mid \alpha \text{ is symplectic } \}. \]
Let $\Sp_{g}(\Z/N)$ be constant group scheme associated to the group of symplectic automorphisms of $(\Z/N)^g \times (\Z/N)^g$. Then $\textbf{S}(A,a,N)$ is a $\Sp_{g}(\Z/N)$ torsor in the fppf topology. It is therefore represented by a scheme which is finite \'etale and Galois with structure group $\Sp_{g}(\Z/N)$ over $S$. 
\end{Proposition}

\begin{proof}
The cover $A[N] \rightarrow S$ is fppf and $\textbf{S}(A,a,N)(A[N])$ is non-empty. The action of $\Sp_{g}(\Z/N)$ on $T$-valued points is clearly transitive and free. Moreover $\textbf{S}(A,a,N)$ is an fppf sheaf for $A[N]$ and $(\Z/N)$ are fppf sheaves and so symplectic isomorphisms between them form an fppf sheaf.  By (\cite{torsor} p363) it is represented by a finite \'etale scheme over $S$.
\end{proof}

\textbf{Remark:} By construction there is a canonical morphism $\textbf{S}(A,a,N) \rightarrow \mathcal{A}_{g,1,N}$. The induced morphism
\[ \mathcal{A}_{g,1} \otimes \Z[1/N] \rightarrow [\mathcal{A}_{g,1,N}/ \Sp_g(\Z/N)] : (A,a)/S \mapsto (\textbf{S}(A,a,N)/S, \textbf{S}(A,a,N) \rightarrow \mathcal{A}_{g,1,N}) \]
is an isomorphism of stacks as can be seen by using descent, although we do not make explicit use of this fact.

\section{Proof of Theorem \ref{twist}}

In this section all schemes will be over $\Z[1/2]$. We observe the following

\begin{Proposition} \label{reg}
Let $R_1 \subset R_2$ be a finite extension of local Noetherian regular rings. Then $R_2$ is free as an $R_1$ module.
\end{Proposition}
\begin{proof}
As $R_1$ is regular, the global dimension of $R_1$ is equal to the dimension of $R_1$ (\cite{serreLoc} p77 Corollary 1).
The projective dimension of $R_2$ is therefore equal to the dimension of $R_1$ minus the depth of $R_2$ over $R_1$ (\emph{ibid.} p75 Proposition 21).
Now $R_2$ and $R_1$ are Cohen-Macaulay as they are regular (\emph{ibid.} p77 Corollary 3) and therefore the depth of $R_2$ over $R_1$ is equal to the dimension of $R_1$ (\emph{ibid.} p63 Proposition 12). Thus the projective dimension of $R_2$ over $R_1$ is zero and $R_2$ is free.
\end{proof}

Let $N \geq 3$ be an integer. The Torelli morphism
\[ \mathfrak{t}_{N}: \curl{M}_{g,N} \rightarrow \curl{A}_{g,1,N} : (C, \alpha) \mapsto (\Jac(C), \lambda_\Theta, \alpha) \]
admits an involution
\[ \tau \map{\curl{M}_{g,N}}{\curl{M}_{g,N}} : (C, \alpha) \mapsto (C, -\alpha). \]
Let $k$ be a field, the Torelli theorem \cite{Weil} states that $C/k$ admits an automorphism mapping to the $[-1]$ automorphism of $\Jac(C)$ if and only if $C$ is hyperelliptic, i.e.
\[  (C, \alpha) \cong (C, -\alpha) \]
if and only if $C$ is hyperelliptic. That is the fixed locus of $\tau$ is equal to the hyperelliptic locus. Let $V_{g,N}$ be the geometric quotient of $\curl{M}_{g,N}$ by the constant group scheme associated to $\langle \tau \rangle$. We note that $V_{g,N}$ exists and yields a finite surjective morphism $q_{g,N} \map{\mathcal{M}_{g,N}}{ V_{g,N}}$ and a morphism $\iota_N \map{V_{g,N}}{\mathcal{A}_{g,1,N}}$ such that $\tor_N = \iota_N \circ q_{g,N}$ (\cite{ab} III. 12). For every field $k$ not of characteristic $2$ the morphism $\iota_{N} \otimes k$ is an immersion (\cite{oortSteen} Theorem 3.1). Thus $\iota$ is unramified; let $Z$ be the closure of its image in $\mathcal{A}_{g,1,N}$. Then let $x \in V_{g,N}$ and let $y \in Z$ be its image and $\mathfrak{p} \in \Spec(\Z[1/N, \zeta_n])$ be the image of $x$ under the structure morphism. The rings $\mathcal{O}_{V_{g,N},x}$ and $\mathcal{O}_{Z,y}$ are local $\Z[1/N, \zeta_n]_{\mathfrak{p}}$ algebras, and as $\iota \otimes k(\mathfrak{p})$ is an immersion and $V_{g,N}$ is integral, by the local criterion for flatness applied to the unique prime $p \in \mathfrak{p} \cap \Z$ the ring $\mathcal{O}_{V_{g,N},x}$ is flat as an $\mathcal{O}_{Z,y}$ module. Thus the morphism $V_{g,N} \rightarrow Z$ is \'etale of degree $1$ and therefore an open immersion, hence $\iota$ is an immersion. In particular given $(A,a)/\Spec(R)$ as in Theorem \ref{twist} the canonical morphism from $\Spec(R)$ to $\mathcal{A}_{g,1,N}$ factors through $V_{g,N}$. We remark that $V_{g,N}$ is Noetherian for an ascending chain of sheaves of ideals on $V_{g,N}$ corresponds to an ascending chain of $\tau$ invariant sheaves of ideals on $\mathcal{M}_{g,N}$, and $\mathcal{M}_{g,N}$ is Noetherian.

\begin{Proposition} \label{free}
Let $H_{g,N}$ denote the image of the hyperelliptic locus in $V_{g,N}$. The sheaf $(q_{g,N})_* \mathcal{O}_{\mathcal{M}_{g,N} \backslash q_{g,N}^{-1}(H_{g,N})}$ is a locally free $\mathcal{O}_{V_{g,N} \backslash H_{g,N}}$-module of rank $2$. The sheaf $(q_{3,N})_* \mathcal{O}_{\mathcal{M}_{3,N}}$ is a locally free $\mathcal{O}_{V_{3,N}}$-module of rank $2$
\end{Proposition}

\begin{proof}
As $V_{g,N}$ is the geometric quotient of $\mathcal{M}_{g,N}$ and $\tau$ is fixed point free away from $q_{g,N}^{-1}(H_{g,N})$, the scheme $V_{g,N} \backslash H_{g,N}$ is Noetherian as remarked, and regular as $\mathcal{M}_{g,N}$ is regular and Noetherian and so the first claim follows from Proposition \ref{reg}. From the classical Torelli theorem we deduce that the map from $V_{3,N}$ to $\mathcal{A}_{3,1,N}^i$ is an isomorphism. The moduli spaces $\curl{M}_{3,N}$ and $\curl{A}_{3,1,N}^i$ are regular and Noetherian and therefore $V_{3,N}$ is regular and Noetherian and the proposition follows from Proposition \ref{reg}.
\end{proof}

\textbf{Remarks: 1)} The second part of Proposition \ref{free} is false for $g > 3$ as then $V_{g,N}$ is not regular and Theorem 13b) p83 of \cite{serreLoc} implies that Proposition \ref{free} is true if and only if $V_{g,N}$ is regular. Indeed Oort and Steenbrink have shown that the tangent space at a hyperelliptic point of $V_{g,N}$ has dimension $g(g+1)/2$. However the Krull dimension of $V_{g,N}$ is clearly $3g-3$.\\
\textbf{2)} We are indebted to Marius van der Put for pointing out Proposition \ref{reg} to us and the fact that it implied Proposition \ref{free} which is the basis of the construction of the twisting functions of Theorem \ref{twist}.  \\

We observe that $\Sp_{g}(\Z/N)$ acts on $\curl{M}_{g,N}, V_{g,N}$ and $\mathcal{A}_{g,1,N}$ and that $\tor_N$ is $\Sp_{g}(\Z/N)$ equivariant. We write $V_{g,N}^\sm$ for the regular locus of $V_{g,N}$. Thus $V_{g,N}^\sm = V_{g,N} \backslash H_{g,N}$ if $g > 3$ and $V_{3,N}^\sm = V_{3,N}$. 

\begin{proof}[Proof of Theorem \ref{twist}]
\textit{Step I} We first construct the twisting functions: \\

Put $S = \Spec(R)$. Fix an integer $N \geq 3$ whose image in $R$ is invertible, for example $N=4$, and set $T_N = \textbf{S}(A,a,N)$ as in Proposition \ref{torsor}. We note that $T_N$ is affine as $T_N/S$ is finite and we write $R_N$ for the affine coordinate ring of $T_N$. The Abelian scheme $(A,a)_{T_N}$ has a universal level $N$ symplectic structure $\alpha$ by definition of $T_N$ and thus there is a morphism
\[ \psi \map{T_N}{V_{g,N}^\sm} \]
corresponding to $(A_{T_N},a_{T_N}, \alpha)$. \\

Put $M_N = T_N \times_{V_{g,N}^\sm} \mathcal{M}_{g,N}$. By construction there is a family of genus $g$ curves $C/M_N$ with a symplectic level $N$ structure $\alpha$ on the Jacobian $\Jac(C)/M_N$ and an isomorphism $((A,a)_{T_N})_{M_N} \cong (\Jac(C),\lambda_\Theta)$ which is $\Sp_{g}(\Z/N)$ equivariant. \\

Then the morphism $M_N \rightarrow T_N$ is affine as $\tor_N$ is affine and the base change of an affine morphism is affine. Thus $M_N$ is affine and we write $R_N^\prime$ for its coordinate ring. Consider an open affine cover $\{ U_i \}$ of $V_{g,N}^\sm$ over which $(q_{g,N})_* \mathcal{O}_{\mathcal{M}_{g,N}}(U_i)$ is free as an $\mathcal{O}_{V_{g,N}^\sm}(U_i)$ module. Consider also an open affine cover $\{ W_{ij} \}$ of $\psi^{-1}(U_i)$. Then the preimages $X_{ij}$ of the $W_{ij}$ in $M_N$ are affine and cover $M_N$ and moreover
\[ \mathcal{O}_{M_N}(X_{ij}) =  \mathcal{O}_{T_N}(W_{ij}) \otimes_{\mathcal{O}_{V_{g,N}^\sm}(U_i)} (q_{g,N})_* \mathcal{O}_{\mathcal{M}_{g,N}}(U_i). \]
Therefore $R_N^\prime$ is locally free of rank $2$ as an $R_N$ module. Furthermore the involution $\tau$ acts on $R_N^\prime$ so that $R_N = (R_N^\prime)^\tau$ and the inclusion $R_N \subset R_N^\prime$ is $\Sp_{g}(\Z/N)$ equivariant. By definition of $T_N$ we have $R = R_N^{\Sp_{g}(\Z/N)}$. Set
\[ R^\prime = (R_{N}^\prime)^{\Sp_{g}(\Z/N)}. \]
Then the map
\[ R^\prime \otimes_R R_N \rightarrow R_N^\prime : a \otimes b \mapsto ab, \]
is an isomorphism, as the category of $R_N$ modules with an $\Sp_{g}(\Z/N)$ action is equivalent to the category of $R$ modules and the functor which takes invariants is inverse to the functor $\otimes_R R_N$ (\cite{ab} III. 12). In particular $M_N/\Spec(R^\prime)$ is \'etale as the base change of an \'etale morphism is \'etale. \\

Now $R^\prime$ is a locally free $R$ module of rank $2$ (\emph{ibid.}). Without loss, we may assume that $R^\prime$ is free as an $R$ module. Then there exists a function $F \in R^\prime$ such that
\[ R^\prime = R \cdot \langle 1, F \rangle \]
and
\[ F^2 = aF + b \]
for some $a,b \in R$. Put
\[ \sqrt{\Delta} = F - a/2. \]
Then $\Delta = F^2 - aF + a^2/4 \in R$. Moreover as $R^\prime_N = R^\prime \otimes_R R_N$ we have
\[ R^\prime_N = R_N[ \sqrt{\Delta}] \] 
and
\[ \tau(\sqrt{\Delta}) = - \sqrt{\Delta}. \]

\textit{Step II} We now show that $\Delta$ is a twisting function for $(A,a)_{T_N}$. \\

Let $\pr_\psi$ be the morphism from $M_N$ to $\mathcal{M}_{g,N}$ induced by $\psi$. If $g =3$ and $x \in M_N(k)$ is hyperelliptic; then $\tau(x) = x$ and therefore $\sqrt{\Delta}(x) = 0$. Likewise for any $x \in M_N(k)$ if $\sqrt{\Delta}(x) = 0$ then $\tau(x) = x$ and therefore $x \in M_N(k)$ is hyperelliptic. \\

Now given $x \in T_N(k)$, if $\Delta(x)$ is a square then $M_N/T_N$ is split above $x$ and there exists a $y \in M_N(k)$ lying above $x$.  Then $\tor_N(\pr_\psi \circ y) = \psi \circ x$ and from the moduli interpretation $(\Jac(C), \lambda_{\Theta})_y = ((A,a)_{T_N})_x$. \\

Let $x \in T_N(k)$ and assume there is a curve $C_0/k$ such that $((A,a)_{T_N})_x = (\Jac(C_0), \lambda_{\Theta_0})$. This means that $\Jac(C_0)$ admits a level $N$ structure $\alpha_0$ and thus
\[ \tor_N^{-1}(\psi \circ x)(k) = \{ (C, \alpha), (C, -\alpha) \} \]
and therefore $\tor_N^{-1}(\psi \circ x)$ and hence $M_N \times_{T_N} x$ is split over $k$. This implies that $\Delta(x)$ is a square. \\

\textit{Step III} We now show that $\Delta$ is a twisting function for $(A,a)/S$. \\

Given $\sigma_1, \sigma_2 \in \Sp_{g}(\Z/N)$ there is a canonical isomorphism $\theta_{\sigma_1, \sigma_2} : \sigma_1^*C \cong \sigma_2^*C$ given by pulling back the corresponding isomorphism on the universal curve over $\mathcal{M}_{g,N}$. Thus the pullbacks to $M_N \times_{\Spec(R^\prime)} M_N$ of $\sigma_1^*C$ and $\sigma_2^*C$ are isomorphic via $\theta_{\sigma_1, \sigma_2,M_N}$ and moreover
\[ \theta_{\sigma_1, \sigma_3, M_N \times M_N} = \theta_{\sigma_2, \sigma_3, M_N \times M_N} \circ \theta_{\sigma_1, \sigma_2, M_N \times M_N}. \]
In other words the $\sigma^*C/M_N$ and the isomorphisms $\theta_{\sigma_1, \sigma_2}$ form descent datum with respect to the \'etale cover $M_N/\Spec(R^\prime)$. As $\mathcal{M}_g$ is a stack, descent datum are effective and there exists a unique curve (up to isomorphism) $\tilde{C}/\Spec(R^\prime)$ whose pull back to $M_N$ is isomorphic with $C/M_N$. Similarly as $\mathcal{A}_{g,1}$ is a stack the isomorphism $((A,a)_{T_N})_{M_N} \cong (\Jac(C),\lambda_\Theta)$ descends to an isomorphism $(A,a)_{\Spec(R^\prime)} \cong (\Jac(\tilde{C}), \lambda_{\tilde{\Theta}})$. \\

Now given $x \in S(k)$ if $\Delta(x)$ is a square then there is a point $y \in \Spec(R^\prime)(k)$ lying above $x$ and $(\Jac(\tilde{C}), \lambda_{\tilde{\Theta}})_y = (A,a)_x$. \\

Assume there is a curve $C_0/k$ such that $(\Jac(C_0), \lambda_{\Theta_0}) \cong (A,a)_x$. As $T_N/S$ is \'etale, the scheme $(T_N)_x$ is isomorphic to the spectrum of a direct sum of fields; as $T_N/S$ is Galois these fields are all isomorphic to a given field $k_N$. That is for some positive integer $s$ we have
\[ (T_N)_x = \Spec( \oplus^s_{i=1} k_N ). \]
By Step II $\Delta$ is a twisting function for $T_N$ and therefore $\Delta(x)$ is a square in $k_N$. Therefore the scheme $(M_N)_{\Spec(R^\prime)_x}$ has exactly two connected components, each of which is isomorphic with $(T_N)_x$. In other words
\[ (M_N)_{\Spec(R^\prime)_x} = \Spec( \oplus_{j=1}^2 \oplus^s_{i=1} k_N ). \]
The scheme $\Spec(R^\prime)_x$ is isomorphic to the spectrum of $k[X]/(X^2 - \Delta(x))$. We therefore have an inclusion of rings
\[ \xymatrix{ \oplus_{j=1}^2 \oplus^s_{i=1} k_N & \ar[l] \oplus^s_{i=1} k_N \\ k[X]/(X^2 - \Delta(x)) \ar[u] & \ar[l] k \ar[u].} \]
Now both vertical arrows are Galois extensions of \'etale $k$-algebras with Galois group $\Sp_{g}(\Z/N)$. Therefore given $y_N \in (M_N)_{\Spec(R^\prime)_x}$ and $x_N \in (T_N)_x$ and $y \in \Spec(R^\prime)_x$ such that $y_N$ lies above $y$, the decomposition group of $y_N/y$ is equal to the decomposition group of $x_N$. The residue field of $y_N$ and $x_N$ is $k_N$ and the residue field of $y$ is $k(\sqrt{\Delta}(x))$. Thus the Galois groups of $k_N/k$ and $k_N/k(\sqrt{\Delta}(x))$ are equal, which implies that $\Delta(x) \in k^2$. \\

Now assume $g =3$ and assume $k = k$ is a field, if $y_N \in M_N(k^s)$ lies above $y \in \Spec(R^\prime)(k)$ then $C_{y_N} \otimes k^s \cong \tilde{C}_y \otimes k^s$. Moreover $\Delta(x)$ is zero if and only if $\sqrt{\Delta}(y_N)$ is zero. Thus the locus of vanishing of $\Delta$ is the hyperelliptic locus on $S$. 
\end{proof}

\section{Line bundles on stacks and Modular forms} \label{mod}
	\begin{Definition}[Mumford \cite{mumPic} p64]
	 	A line bundle $L$ on a stack $\curl{X}$ is given by the following data:
	\begin{itemize}
	 \item for each section $x \map{S}{\curl{X}}$ an $\mathcal{O}_S$ module $L(x)$ on $S$
	\item for each morphism $f \map{S^\prime}{S}$ an isomorphism $\psi_{f,x} \map{f^*L(x)}{L(x \circ f)}$.
	\end{itemize}
	subject to the compatibility condition
	\[ \psi_{f \circ g,x} \circ \xi =  \psi_{g, x \circ f} \circ g^* \psi_{f,x} \]
	where $\xi$ is the natural isomorphism between $g^*f^*L( x)$ and $(f \circ g)^*L(x)$. \\
	A global section of a line bundle $L$ on a stack $\curl{X}$ is given by the following data: for each section $x \map{S}{\curl{X}}$ there is a section $s_x \in \cH^0(S,L(x))$, subject to the condition $\psi_{f,x}(f^*s_x) = s_{x \circ f}$. 
	\end{Definition}
	
	For example $\mathcal{O}_{\mathcal{X}}$ is given by $\mathcal{O}_{\mathcal{X}}(x \map{S}{\mathcal{X}}) = \mathcal{O}_S$. We write $\cH^0(\mathcal{X}, L)$ for the $\cH^0 (\mathcal{X}, \mathcal{O}_{\mathcal{X}})$ module of global sections of $L$. \\
	
	Given global sections $s_0$ and $s_1$ of a line bundle $L$ on a stack $\mathcal{X}$ if for each section $x \map{S}{\curl{X}}$ the Cartier divisor of zeroes of $s_0$ is equal to the Cartier divisor of zeroes of $s_1$ then there is a function $\lambda_x \in \mathcal{O}_S^*(S)$ such that
	\[ \lambda_x s_0 = s_1. \]
	In particular the $\lambda_x$ are unique as $L(x)$ is a line bundle and therefore the $\lambda_x$ define a global section $\lambda$ of the structure sheaf $\mathcal{O}_{\curl{X}}$ of $\curl{X}$ which is a unit. \\
	
	Given a scheme $X$ and a morphism $f \map{\mathcal{X}}{X}$ which is initial in the category of maps from $\mathcal{X}$ to schemes we say that $X$ is a coarse moduli space for $\mathcal{X}$ if the isomorphism classes of $\mathcal{X}(k)$ are in bijection via $f$ with $X(k)$ for every algebraically closed field $k$. Then the natural injection
	\[ \cH^0(X, \mathcal{O}_X) \subset \cH^0(\mathcal{X}, \mathcal{O}_{\mathcal{X}}) \]
 is an equality: an element of the right hand side corresponds to a morphism $f \map{\mathcal{X}}{\Aff^1}$ which factors uniquely through $X$. The stacks $\mathcal{A}_{g,1}/\Z$ and $\mathcal{M}_{g}/\Z$ both have coarse moduli spaces (\cite{GIT} Theorem 7.10, Corollary 7.14), which are Zariski locally quasi-projective and therefore separated and of finite type over $\Z$. \\

Let $f \map{X}{S}$ be either an Abelian scheme or a family of genus $g$ curves.  For lack of reference we prove
 \begin{Proposition} \label{free2}
 Let $f \map{X}{S}$ be as above, then the sheaf $f_*\Omega_{X/S}$ is locally free.
 \end{Proposition}
 \begin{proof}
 In case $X/S$ is an Abelian scheme the natural morphism $f_*\Omega_{X/S} \rightarrow O^*\Omega_{X/S}$ is an isomorphism and the Proposition follows. In case $X/S$ is a family of curves (i.e. smooth proper with geometrically connected fibres of dimension $1$), then $\mathcal{O}_X$ is flat over $\mathcal{O}_S$. Thus higher direct images of $\mathcal{O}_X$ under $f$ commute with arbitrary base change (\cite{EGAIII2} 7.7.5.3 for $Y=S$ and $\mathcal{P}^\bullet$ equal to the complex supported in degree $0$ determined by $\mathcal{O}_X$). Therefore $R^1f_*\mathcal{O}_C$ is locally free (\emph{ibid} apply Propositions 7.7.10a) \& 7.8.4 to Proposition 7.8.5). By relative Serre duality $f_*\Omega_{X/S}$ is locally free (\cite{liu}).
 \end{proof} 
 
 The Hodge bundle of $f$ is the sheaf
	\[ \omega_{X/S} := \det(f_*\Omega_{X/S}). \]
	Let $\rho$ be a positive integer; a section $s$ of $\omega_{X/S}^\rho$ is called a modular form of weight $\rho$. We have line bundles $\omega_{A/\curl{A}_{g,1}}^\rho$ and $\omega_{C/ \curl{M}_{g}}^\rho$ on $\curl{A}_{g,1}$ and $\curl{M}_g$ respectively defined by the corresponding sheaves of modular forms because of basic properties of differentials under base change. Sections of these line bundles are called Katz-Siegel and Katz-Teichm\"uller modular forms respectively (cf. \cite{Katz}). \\
	
\begin{Proposition}
If $f \map{C}{S}$ is a family of genus $g$ curves and $\pi \map{\Jac(C)}{S}$ is its Jacobian then $\omega_{C/S}$ is isomorphic with $\omega_{\Jac(C)/S}$.
\end{Proposition}
\begin{proof}
As $\Jac(C)/S$ is a group scheme we have an isomorphism
\[ \pi_*\Omega_{\Jac(C)/S} = O^*\Omega_{\Jac(C)/S}. \]
Interpreting sections of $R^1f_*\curl{O}_C$ as Cech $1$-cocycles shows that 
\[ R^1f_*\curl{O}_C (U) = \Hom( O^*\Omega_{\Jac(C)/S}(U), \mathcal{O}_{S}(U)). \]
Then Serre duality for smooth and proper morphisms shows that (\cite{liu})  
\[ \pi_*\Omega_{\Jac(C)/S} \cong f_*\Omega_{C/S}. \]
\end{proof}

Thus the pullback of a Katz-Siegel modular form via the Torelli morphism is a Katz-Teichm\"uller modular form.

\section{The modular form $\Disc$} \label{disc}

\begin{Definition}
The classical family of plane quartics is the family of quartics $f \map{Q}{\Pro^{14}}$ where $Q$ is the subscheme of $\Pro^2_{\Pro^{14}}$
  defined by the equation
  \[ ( F := \sum_{i+j+l = 4} a_{ijl} X^i Y^j Z^{l} ) = 0. \]
\end{Definition}

 \begin{Definition}[Classical discriminant]
  Let $F_X$, $F_Y$ and $F_Z$ be the derivatives of $F$ with respect to $X, Y$ and $Z$.
  The classical discriminant of $f \map{Q}{\Pro^{14}}$ is defined to be the resultant of $F_X,
  F_Y$ and $F_Z$; it is a section of $\omega_{Q/\Pro^{14}}^9$.
\end{Definition}

  We recall the main properties of the resultant (see \cite{Lang} pp388-404 for details). Let $\{ U_j \}$ be the standard affine cover of $\Pro^{14}$. The locus of $(a_{ijl}) \in U_i$ such that $F_X,F_Y$ and $F_Z$ have a common zero is defined by the single homogenous equation $\Res_{U_i}(F_X,F_Y,F_Z)$ which generates a prime ideal. Moreover there is a canonical choice of $\Res_{U_i}(F_X,F_Y,F_Z)$. Finally if $B$ is a $3$ by $3$ matrix with entries in $U_i$, then 
\[ \Res(F_X~\circ~B^{-1},~F_Y~\circ~B^{-1},~F_Z~\circ B^{-1}) = \det(B)^9 \Res(F_X, F_Y,F_Z). \]  
  Therefore the $\Res_{U_i}(F_X,F_Y,F_Z)$ glue to form a unique global section $\Res(F_X,F_Y,F_Z)$ of $\omega_{Q/\Pro^{14}}^9$.

Let $f \map{C}{S}$ be any family of genus $3$ curves. There is a affine open cover $\{ V_i \}$ of $S$ such that $f_*\Omega_{C_{V_i}/V_i}$ and $f_*\Sym^4(\Omega_{C_{V_i}/V_i})$ are free (by analogous reasoning to Proposition \ref{free2}). Put $V_i = \Spec(R_i)$ and let $A_i$ be the graded $R_i$ algebra given by the direct sum of the tensor powers of the $R_i$ modules $f_*\Omega_{C_{V_i}/V_i}(V_i)$. We then have a morphism
\[ j_i \map{C_{V_i}}{\Proj(A_i)} \]
and
\[ \Proj(A_i) \cong \Pro^2_{V_i}. \]
Let $I_i$ be the homogenous ideal corresponding to the image of $C_{V_i}$. The degree $4$ elements $I_i^{(4)}$ of $I_i$ are equal to the kernel of the surjective $R_i$ linear morphism
\[  \Sym^4(f_*\Omega_{C_{V_i}/V_i})(V_i) \rightarrow f_*\Sym^4(\Omega_{C_{V_i}/V_i})(V_i). \]
Now cohomology and base change commute for $\Sym^4(\Omega_{C_{V_i}/V_i})$ thus for each $s \in V_i(k)$
\[ \dim_{k} \cH^0(C_s, \Sym^4(\Omega_{C_s/k})) = \dim_k f_*\Sym^4(\Omega_{C_{V_i}/V_i})_s. \]
Thus by Riemann-Roch for curves over fields
\[ \mathrm{rank}_{R_i}(I_i^{(4)}) = 15 -  \dim_{k} \cH^0(C_s, \Sym^4(\Omega_{C_s/k})) = 1. \]
Thus there is a unique degree $4$ element $F_i \in I_i$ up to $R_i^*$-scaling. There is a morphism $\iota_i : V_i \rightarrow \Pro^{14}$ such that $F_i$ the pullback of $F$ with $\iota_i$ (up to scaling by an element of $R_i^*$). The $\iota_i^*\Res(F_X,F_Y,F_Z)$ are independent of $\iota_i$ and therefore glue to form a global section $\Disc_{C/S}$ of $\omega_{C/S}^9$. Moreover by construction given a morphism $h \map{S^\prime}{S}$ we have a canonical isomorphism $\psi_h \map{h^*\omega_{C/S}}{\omega_{C_{S^\prime}/S^\prime}}$ and $\Disc_{C_{S^\prime}/S^\prime} = \psi_h(h^*\Disc_{C/S})$.

\begin{Definition}
The discriminant modular form $\Disc$ is the Katz-Teichm\"uller modular given by the $\Disc_{C/S}$ as above.
\end{Definition}

\begin{Proposition}
 	$\Disc$ vanishes on the hyperelliptic locus.
\end{Proposition}
\begin{proof}
If $C/k$ is hyperelliptic, then the hyperelliptic structure morphism $f \map{C}{\Pro^1}$ defines an element of the function field $k(C)$ whose polar divisor is equal to half the canonical class. Therefore we have a basis for $\Omega_{C/k}$ given by
\[ \Omega_{C/k} = k \langle 1, f, f^2 \rangle. \]
Thus the corresponding quartic for $C/k$ is the square of a conic and is therefore singular.
 \end{proof}

\section{The modular form $\chi_{18}$}

Let $\mathfrak{h}_3$ be the Siegel upper half space of degree $3$: that is the complex domain of $3 \times 3$ symmetric complex matrices with positive definite imaginary part. Let $\pi^\an \map{A^\an}{\mathfrak{h}_3}$
be the universal complex analytic Abelian threefold over $\mathfrak{h}_3$. Analytic uniformisation provides a family $f \map{\C^3 \times \mathfrak{h}_{3}}{\mathfrak{h}_3}$ of complex vector spaces, and a family $h \map{\Lambda}{\mathfrak{h}_3}$ of lattices of rank $6$ whose fibre
above a point $\tau$ is
\[ h^{-1}(\tau) := \Z^3 + \tau \Z^3 \subset f^{-1}(\tau) = \C^3 \times \{ \tau \} \]
so that
\[ A^\an \cong (\C^3 \times \mathfrak{h}_{3})/ \Lambda \]
as complex analytic manifolds above $\mathfrak{h}_3$. \\

Let $z = (z_1(\tau), z_2(\tau) , z_3(\tau))$ denote the standard coordinates of $\C^3 \times \{ \tau \} \subset \C^3 \times \mathfrak{h}_{3}$. The Riemann theta function is the holomorphic function on $\C^3 \times \mathfrak{h}_3$ defined by the formula
\[ \vartheta(z; \tau) := \sum_{n \in \Z^3} \mathrm{e}^{ \pi  \sqrt{-1} n^t \tau n + 2 \pi \sqrt{-1} n \cdot z}. \]
We have a canonical level $2$ structure on $\pi^\an \map{A^\an}{\mathfrak{h}_3}$ given by the standard coordinates
on $\Z^3 + \tau \Z^3$. Using this we identify a $2$-torsion point $p \in A^\an[2](\mathfrak{h}_3)$ with a unique element $m+\tau m^\prime
\in 1/2 (\Z^3 + \tau \Z^3)$ such that $m$ and $m^\prime$ have entries equal to either $0$ or $1/2$. The analytic theta nulls are the holomorphic functions on $\mathfrak{h}_3$ given by the formulas
\[ \vartheta \pair{m}{m^\prime}(0; \tau) := \vartheta(m + \tau m^\prime; \tau). \]
The function $p \mapsto 4m \cdot m^\prime \mod 2$ defines a quadratic form with values in $\Z/2 \Z$. Thus an analytic theta null
\[ \vartheta \pair{m}{m^\prime}(0; \tau) \]
is called even if $4m \cdot m^\prime$ is even. 
\begin{Definition} Put
\[ \chi_{18}^\hol(\tau) := \prod_{4 m \cdot m^\prime \in 2 \Z} \vartheta \pair{m}{m^\prime}(0; \tau),  \]
and
\[ \chi_{18} = \chi_{18}^\hol (2 \pi )^{54} (dz_1 \wedge dz_2 \wedge dz_3)^{18}. \]
\end{Definition}

Igusa has shown (\cite{igusa} Lemmas 10 and 11) that $\chi_{18}^\hol$ transforms as a weight $18$ Siegel modular form under the group $\PSp_3(\Z)$, and vanishes on the hyperelliptic locus. By the Koecher principle, this means that $\chi_{18}$ is a modular form on a smooth toroidal compactification of $\mathcal{A}_{3,1} \otimes \C$ and therefore an algebraic modular form by GAGA (cf. \cite{FaltChai} p141). Moreover $\chi_{18}$ is a modular form on $\mathcal{A}_{3,1}/ \Z$ by the $q$-expansion principle (\cite{FaltChai} p140). \\

For the reader unhappy with applying GAGA to a stack, despite the authority of Faltings-Chai, we offer an alternative argument. First of all $\chi_{18}$ defines a modular form on the analytic moduli space $\mathcal{A}_{3,1,4}^\an$, which by the Koecher principle, extends to a smooth toroidal compactification of $\mathcal{A}_{3,1,4}^\an$ and is therefore analytic by GAGA. Now if $(A,a)/S$ is any principally polarised Abelian scheme over a $\C$ scheme $S$, the cover $\textbf{S}(A,a,4) \rightarrow S$ (Proposition \ref{torsor}) is Galois with Galois group $\Sp_3(\Z/4)$. Then Igusa's Lemma shows that the pull back of $\chi_{18}$ to $T$ is $\Sp_3(\Z/4)$ invariant and so descends to $S$. Thus $\chi_{18}$ defines an algebraic modular form on $\mathcal{A}_{3,1} \otimes \C$.

\section{Proof of Theorem \ref{chidisc}}
\emph{Step I:}
The global regular functions of $\mathcal{M}_3$ are equal to the global regular functions of its coarse moduli space $M_3$ (see Section 4). Moreover $M_3 \otimes \Q$ is isomorphic with the coarse moduli space of indecomposable principally polarised Abelian $3$-folds $A_{3,1}^i \otimes \Q$ (\cite{oortSteen} Theorem 3.2). The coarse moduli space $A_{3,1}^i \otimes \Q$ is quasi-projective with boundary of codimension $2$ under the Satake embedding. Therefore using the theorem of Bertini (\cite{jou} p89, Th\'eor\`eme 6.10), any two closed points on $A_{3,1}^i \otimes \Q$ may be connected by a complete curve $C$ lying within $A_{3,1}^i \otimes \Q$, and therefore $A_{3,1}^i \otimes \Q$ has only the field of constants $\Q$ as global regular functions. As $\Z \subset \Q$ is flat, and $M_3$ is separated and of finite type over $\Z$, we have (\cite{egaIII} 1.4.15) that $\cH^0(M_3,  \mathcal{O}_{M_{3}}) = \Z$. \\

\emph{Step II:}
The Katz-Teichm\"uller modular forms $\tor^*\chi_{18}$ and $\Disc^2$ both have weight $18$. Assume there is a section $x \map{S}{\mathcal{M}_3}$ such that for all $\gamma \in \Gamma(S, \mathcal{O}_S)$ we have $(\tor^*\chi_{18})_x \neq \gamma \Disc^2_x$. Then there exists an open affine $\Spec(R) \subset S$ such that for all $\gamma \in R$ we have $(\tor^*\chi_{18})_{x,R} \neq \gamma \Disc^2_{x,R}$. Moreover, since for any two distinct primes $p,l \in \Z$ the open affines $\Spec(R[1/l])$ and $\Spec(R[1/p])$ cover $\Spec(R)$, we may without loss assume that we have a rational prime $p \geq 3$ such that $1/p \in R$. Now let $C/\Spec(R)$ be the family of genus three curves corresponding to $R$. Let $T^\prime = \textbf{S}(\Jac(C), \lambda_\Theta,p)$ as in Proposition \ref{torsor}. Then if the pullbacks of $(\tor^*\chi_{18})_{x,R}$ of $\Disc^2_{x,R}$ to $T^\prime$ satisfy $(\tor^*\chi_{18})_{T^\prime} = \gamma \Disc_{T^\prime}^2$ for some $\gamma \in \Gamma(T^\prime, \mathcal{O}_{T^\prime})$ then $\gamma$ is invariant under the $\Sp_3(\Z/p)$ and therefore lies in $R$. Therefore for all $\gamma \in \Gamma(T^\prime, \mathcal{O}_{T^\prime})$ we have $(\tor^*\chi_{18})_{T^{\prime}} \neq \gamma \Disc_{T^{\prime}}^2$. Now the Jacobian of $C_{T^\prime}$ admits a level $p$ structure, and therefore there is a morphism $\psi \map{T^\prime}{\mathcal{M}_{3,p}}$ and
\[ (\tor^*\chi_{18})_{T^\prime} = \psi^* (\tor^*\chi_{18})_{\mathcal{M}_{3,p}} \]
and
\[ \Disc_{T^\prime}^2 = \psi^*\Disc_{\mathcal{M}_{3,p}}^2. \]
Let $C_p$ be the universal curve over $\mathcal{M}_{3,p}$. As $\Z[1/p, \zeta_p] \subset \C$ is flat and $\mathcal{M}_{3,p}$ is separated and of finite type over $\Z[1/p, \zeta_p]$ we have (\emph{ibid.}) that
\[ \cH^0( \mathcal{M}_{3,p}, \omega_{C_p/\mathcal{M}_{3,p}} ) \otimes \C = \cH^0( \mathcal{M}_{3,p} \otimes \C, \omega_{C_{p,\C}/\mathcal{M}_{3,p} \otimes \C} ). \]
Now $(\tor^*\chi_{18})_{\C}$ and $\Disc^2_{\C}$ have the same divisor of zeroes on $\mathcal{M}_{3} \otimes \C$, namely the hyperelliptic locus. Thus there is a $\gamma \in \C$ such that
\[ (\tor^*\chi_{18})_{\mathcal{M}_{3,p}} \otimes 1 =\Disc_{\mathcal{M}_{3,p}}^2 \otimes \gamma. \]
Therefore $\gamma \in \Z[1/p,\zeta_p]^*$, and
\[ (\tor^*\chi_{18})_{\mathcal{M}_{3,p}}  = \gamma \Disc_{\mathcal{M}_{3,p}}^2, \]
which is a contradiction. Thus for all $x \map{S}{\mathcal{M}_3}$ there exists a $\gamma_x \in \Gamma(S, \mathcal{O}_S)$ such that
\[ (\tor^*\chi_{18})_x  = \gamma_x \Disc_{x}^2.  \]
Thus the $\gamma_x$ define an element $\gamma \in \cH^0(M_3,  \mathcal{O}_{M_{3}})^*$ such that
\[ (\tor^*\chi_{18})  = \gamma \Disc_{\mathcal{M}_{3,p}}^2. \]
Therefore
\[ \gamma \in \cH^0(M_3,  \mathcal{O}_{M_{3}})^* = \Z^*. \]

\section{Proof of Theorem \ref{twistchi}}
Let $R$ be a ring containing $1/2$ equipped with an inclusion $R \subset \C$. Let $(A,a)/\Spec(R)$ be an Abelian scheme of relative dimension $3$ with an indecomposable principal polarisation. Let $\Delta$ be a twisting function for $(A,a)$. Assume that $\Omega_{A/R}$ free of rank $3$ and let $\{ \xi_1, \xi_2, \xi_3 \}$ be a basis. Let $s \map{\Spec(R)}{\mathcal{A}_{3,1}^i}$ be the section corresponding to $(A,a)$. Then
\[ \Delta (\xi_1 \wedge \xi_2 \wedge \xi_3)^{18} \]
and
\[ (\chi_{18})_s \]
are both global sections of $\omega_{A/\mathcal{A}_{3,1}^i}(s)$ with the same divisor of zeroes. Therefore there is a $\kappa \in R^*$ such that 
\begin{equation*}
 \kappa \Delta (\xi_1 \wedge \xi_2 \wedge \xi_3)^{18} = (2\pi \sqrt{-1})^{54} \chi_{18}^{\hol} (dz_1 \wedge dz_2 \wedge dz_3)^{18}.
\end{equation*}

Let $B(A,a)$ denote the set of isomorphisms between $\cH_1((A_\C)^\an, \Z)$ and $\Z^6$. Let $M_3(\C)$ be the set of $3 \times 3$ matrices with values in $\C$. We define functions
\[ \Omega_1 \map{B(A,a)}{M_3(\C)} \]
and
\[ \Omega_2 \map{B(A,a)}{M_3(\C)} \]
as follows. For an element $\eta \in B(A,a)$ we let $\eta_i = \eta^{-1}(e_i)$ where $e_i$ is the $i$th standard basis element of $\Z^6$. Then the entries of $\Omega_1(\eta)$ are obtained by integrating the $\xi_i$ along the $\eta_j$ with $j \in \{1,2,3 \}$. Likewise the entries of $\Omega_2(\eta)$ are obtained by integrating the $\xi_i$ along the $\eta_j$ with $j \in \{4,5,6 \}$. \\

Set
\[ \tau(\eta) = \Omega_1(\eta) \cdot \Omega_2(\eta)^{-1}. \]
Then $\tau$ is a bijection between $B(A,a)$ and the elements of $\mathfrak{h}_3$ corresponding to $(A_\C,a_\C)^\an$. \\

Integration defines the non-degenerate Poincar\'e pairing
\[ \cH^1_{dR}((A_\C)^\an, \R) \otimes (\cH_1((A_\C)^\an, \Z) \otimes \R) \rightarrow \R \]
and we let $dz_i^\eta$ denote the differential dual to $\eta_i \otimes 1$. \\

By construction we have the identity
\begin{equation*}
 \left( \begin{array}{c} \xi_1 \\ \xi_2 \\ \xi_3 \end{array} \right) = \Omega_1(\eta) \left( \begin{array}{c} dz_1^\eta \\ dz_2^\eta \\ dz_3^\eta \end{array} \right). 
\end{equation*}

Let $x$ be the image of $0$ under the morphism $\Spec(\C) \rightarrow \Spec(R)$; we note that $x$ is the generic point of $\Spec(R)$, as $\Spec(R)$ is integral and therefore irreducible. The value of $\Delta$ at $x$ is equal to $\Delta$ as $R \subset R_x$. Therefore $\Delta$ is given by evaluating 
\begin{equation*}
 (2\pi \sqrt{-1})^{54} \chi_{18}^{\hol} (dz_1 \wedge dz_2 \wedge dz_3)^{18} (\xi_1 \wedge \xi_2 \wedge \xi_3)^{-18}
\end{equation*}
at a choice of $\eta \in B(A,a)$. Thus
\begin{eqnarray*}
 \kappa \Delta & = & (2\pi \sqrt{-1})^{54} \chi_{18}^{\hol} (dz_1 \wedge dz_2 \wedge dz_3)^{18} (\xi_1 \wedge \xi_2 \wedge \xi_3)^{-18}(\eta) \\
 & = &  (2\pi \sqrt{-1})^{54} \chi_{18}^{\hol}(\tau(\eta)) (dz_1^\eta \wedge dz_2^\eta \wedge dz_3^\eta)^{18} (\xi_1 \wedge \xi_2 \wedge \xi_3)^{-18} \\
 & = &  \frac{(2\pi \sqrt{-1})^{54} \chi_{18}^{\hol}(\tau(\eta)) }{\det(\Omega_1(\eta))^{18}}
\end{eqnarray*} 

We note the quantity on the last line is invariant under symplectic transformations, and it depends on $A/R$, as $\Omega_1(\eta)$ depends on $A/R$. If $A^\tw$ is the $[-1]$ twist of $A$ over $R[\sqrt{D}]$ with $D \in R \backslash R^2$ then the global differentials of $A^\tw$ in $\Omega_{A \otimes R[\sqrt{D}]/R[\sqrt{D}]}(A \otimes R[\sqrt{D}])$ are invariant under the $R$-linear automorphism given by
\[ \xi_i \mapsto - \xi_i \]
and
\[ \sqrt{D} \mapsto - \sqrt{D}. \]
Therefore the $\sqrt{D} \xi_i \in \Omega_{A \otimes R[\sqrt{D}]/R[\sqrt{D}]}(A \otimes R[\sqrt{D}])$ form an $R$ basis for the differentials on $A^\tw$. Thus $\det(\Omega_1(\eta))^{18}_A = D^{-27} \det(\Omega_1(\eta))^{18}_{A^\tw}$. \\

It remains to show that $\kappa/ \gamma$ is a square in $k$ for each $x \in \Spec(R)(k)$. If $\Delta(x)$ is a square, then there exists a curve $C/k$ in the fibre of $\tor^{-1}((A,a)_x)$ and thus
\[ \kappa/\gamma \Delta(x) (\xi_1 \wedge \xi_2 \wedge \xi_3)^{18}_x = \Disc_{C/k}^2. \]
Thus if $\Delta(x)$ is a square then $\kappa/\gamma \in k^2$. If $\Delta(x)$ is not a square in $k$ the $[-1]$ twist of $(A,a)_x$ over $k(\sqrt{\Delta}(x))$ is a Jacobian and thus as elements of $k$
\[ \Delta^{-27} \frac{(2\pi \sqrt{-1})^{54} \chi_{18}^{\hol}(\tau(\eta)) }{\det(\Omega_1(\eta))^{18}} = \gamma \Disc_{C/k}^2. \]
And therefore
\[ \kappa/\gamma \Delta(x) (\xi_1 \wedge \xi_2 \wedge \xi_3)^{18}_x = \Delta^{-27}(x) \Disc_{C/k}^2 \]
and therefore $\kappa/ \gamma$ is a square in $k$.

\section{Alternative Proof of Theorem \ref{twist} for genus $3$}
Let $R$ be a $\Z[1/2]$ algebra. Let $(A,a)/\Spec(R)$ be an indecomposable principally polarised Abelian scheme of relative dimension $3$. Assume that $\Omega_{A/\Spec(R)}$ is free and choose an isomorphism
\[ \phi \map{\omega_{A/R}(R)}{R}. \]
Fix a basis $\xi_i$ of $\Omega_{A/R}(A)$. We claim the function
\[ \Delta = \gamma^{-1} \phi(\chi_{18}) \phi(\xi_1 \wedge \xi_2 \wedge \xi_3)^{-18} \]
is a twisting function. \\

In the first place $\chi_{18}$ vanishes on the hyperelliptic locus by Theorem \ref{chidisc}. \\

Given $x \in \Spec(R)(k)$, if $(A,a)_x$ is a Jacobian, then
\[ \Delta(x) = \phi(\Disc^2_x) \phi(\xi_1 \wedge \xi_2 \wedge \xi_3)^{-18}, \]
is a square. On the other hand, if $(A,a)_x$ is not a Jacobian, there is a $D \in k \backslash k^2$ such that the $[-1]$ twist $(A,a)_\epsilon$ of $(A,a)_x$ over $k[\sqrt{D}]$ is a Jacobian. Then $\{ \sqrt{D} \xi_i \mid i \in \{ 1,2,3 \} \}$ is an $k$-basis of differentials for $(A,a)_\epsilon$ in $\Omega_{(A,a) \otimes k[\sqrt{D}]}(A \otimes k[\sqrt{D}])$.
Therefore
\[ \Delta(x) = D^{27} \Delta( (A,a)_\epsilon) \]
is not a square.

 \section*{Acknowledgements}
  	We thank Christophe Ritzenthaler and Gilles Lachaud for sharing their work with us at an early stage. We thank Robert Carls, Everett Howe, Ching-Li Chai, David Kohel,
Gerard van der Geer, Bert van Geemen, Jaap Top, Marius van der Put and Jean-Pierre Serre for comments upon earlier versions of this work.

\end{document}